\documentclass{birkjour}
\usepackage{amssymb, amsmath}
\usepackage{enumerate}
\usepackage{color, soul, cancel}

\theoremstyle{plain}

\newtheorem{theorem}{Theorem}

\newtheorem{definition}[theorem]{Definition}

\theoremstyle{remark}

\newtheorem{example}{Example}

\newtheorem{remark}[theorem]{Remark}

\newcommand{\bit}{\begin{itemize}}
\newcommand{\eit}{\end{itemize}}
\newcommand{\ben}{\begin{enumerate}}
\newcommand{\een}{\end{enumerate}}
\newcommand{\be}{\begin{equation}}
\newcommand{\ee}{\end{equation}}
\newcommand{\ba}{\begin{array}}
\newcommand{\ea}{\end{array}}

\newcommand{\abs}[1]{\left|#1\right|}
\newcommand{\norm}[1]{\left|\left|#1\right|\right|}

\newcommand{\dt}{\mathrm{d}t}

\newcommand{\eps}{\varepsilon}

\newcommand{\inner}[2]{\left\langle #1,#2 \right\rangle}

\newcommand\cA{\mathcal A}
\newcommand\cB{\mathcal B}

\newcommand\cH{\mathcal H}
\newcommand\cK{\mathcal K}

\newcommand\C{\mathbb C}
\newcommand\R{\mathbb R}
\newcommand\N{\mathbb N}

\newcommand{\ler}[1]{\left( #1 \right)}

\begin{document}

%
%
%
%
%
%
%
%
%

\title[Centrality and local convexity]
 {Characterizations of centrality by local convexity of certain functions on $C^*$-algebras}

\author[D. Virosztek]{D\'aniel Virosztek}

\address{Functional Analysis Research Group, Bolyai Institute\\
         University of Szeged\\
         H-6720 Szeged, Aradi v\'ertan\'uk tere 1.\\ Hungary}

\email{virosz89@gmail.com}

\thanks{The author was partially supported by the Hungarian National Research, Development and Innovation Office – NKFIH (grant no. K124152).}


\date{October 31, 2017}

\begin{abstract}
We provide a quite large function class which is useful to distinguish central and non-central elements of a $C^*$-algebra in the following sense: for each element $f$ of this function class, a self-adjoint element $a$ of a $C^*$-algebra is central if and only if $f$ is locally convex at $a.$
\end{abstract}

\maketitle
\section{Introduction}

\subsection{Motivation}
Connections between algebraic properties of $C^*$-algebras and some essential properties of functions defined on them by functional calculus have been investigated widely.
\par
The first results concern the relation between the commutativity of a $C^*$-algebra and the monotonicity (with respect to the order induced by positivity) of certain functions defined on the positive cone of it. It was shown by \emph{Ogasawara} in 1955 that a $C^*$-algebra is commutative if and only if the map $a \mapsto a^2$ is monotone increasing on its positive cone \cite{oga}. Later on, \emph{Pedersen} provided a generalization of Ogasawara's result for any power function $a \mapsto a^p$ with $p>1$ \cite{pedersen}. More recently, \emph{Wu} proved that the exponential function is also useful to distinguish commutative and non-commutative $C^*$-algebras in the above sense \cite{wu}, and in 2003, \emph{Ji} and \emph{Tomiyama} described the class of all functions that can be used to decide whether a $C^*$-algebra is commutative or not \cite{ji-tom}.
\par
Some "local" results were also obtained in this topic. First, \emph{Moln\'ar} showed that a self-adjoint element $a$ of a $C^*$-algebra is central if and only if the exponential function is \emph{locally} monotone at $a$ \cite{mol-char}. Later on, we managed to provide a quite large class of functions (containing all the power functions with exponent greater than $1$ and also the exponential function) which has the property that each element of this function class can distinguish central and non-central elements via local monotonicity \cite{dv-mon}.
\par
Investigating the connections between the commutativity of a $C^*\!$-algebra (or locally, the centrality of an element) and the global (or local) convexity property of some functions is of a particular interest, as well.
\par
In 2010, {\it Silvestrov, Osaka and Tomiyama} showed that a $C^*$-algebra $\cA$ is commutative if and only if there exists a convex function $f$ defined on the positive axis which is not convex of order $2$ (that is, it is not convex on the $C^*$-algebra of the $2 \times 2$ matrices) but convex on $\cA$ \cite[Thm. 4.]{proc-eston}.
\par
Motivated by the above mentioned result in \cite{proc-eston}, the main aim of this paper is to provide a large class of functions which have the property that they are locally convex only at central elements, that is, they characterize central elements by local convexity.

\subsection{Basic notions, notation}

Throughout this paper, $C^*$-algebras are always assumed to be unital.
The spectrum of an element $a$ of the $C^*$-algebra $\cA$ is denoted by $\sigma(a).$ The symbol $\cA_s$ stands for the set of all self-adjoint elements of $\cA.$ A self-adjoint element of a $C^*$-algebra is called \emph{positive} if its spectrum is contained in $[0,\infty).$ The order induced by positivity on the self-adjoint elements is defined as follows: $a \leq b$ if $b-a$ is positive. In the sequel, the symbol $\cH$ stands for a complex Hilbert space and $\cB(\cH)$ denotes the algebra of all bounded linear operators on $\cH.$ The inner product on a Hilbert space is denoted by $\inner{\cdot}{\cdot}$ and the induced norm is denoted by $\norm{\cdot}.$ If $u$ and $v$ are elements of a Hilbert space, the symbol $u \otimes v$ stands for the linear map $z \mapsto \inner{z}{v} u.$

\section{The main theorem}

In this section we provide the main result of this paper. In order to do so, we first need a definition.

\begin{definition}[Local convexity] \label{def:loc-conv}
Let $\cA$ be a $C^*$-algebra and let $f$ be a continuous function defined on some open interval $I \subset \R.$ Let $a \in \cA_s$ with $\sigma(a) \subset I.$ We say that $f$ is \emph{locally convex} at the point $a$ if for every $b \in \cA_s$ such that $\sigma(a+b) \cup \sigma(a-b) \subset I$ we have
$$
f(a) \leq \frac{1}{2}\ler{f(a+b)+f(a-b)}.
$$
\end{definition}

\begin{remark} \label{rem:mid-point}
Note that in fact the above definition is the definition of the \emph{mid-point convexity}. However, in this paper every function is assumed to be continuous, so there is no difference between mid-point convexity and convexity.
\end{remark}

Now we are in the position to present the main result of the paper.

\begin{theorem}\label{thm:main}

Let $I \subset \R$ be an open interval and let $f$ be a convex function in $C^2(I)$ such that the second derivative $f''$ is strictly concave on $I.$
Let $\cA$ be a $C^*$-algebra and let $a \in \cA_s$ be such that $\sigma(a) \subset I.$
Then the followings are equivalent.
\ben [(1)]
\item The element $a$ is central, that is, $ab=ba$ for every $b \in \cA.$ \label{tul:cent}
\item The function $f$ is locally convex at $a$. \label{tul:loc-conv}
\een
\end{theorem}

\begin{example} \label{ex: jo-fvek}
On the interval $I=(0, \infty)$ the functions $f(x)=x^p \, (2<p<3)$ satisfy the conditions given in Theorem \ref{thm:main}. That is, these functions are useful to distinguish central and non-central elements via local convexity.
\end{example}

\section{Proof of the main theorem}

This section is devoted to the proof of Theorem \ref{thm:main}. We believe that some of the main ideas of the proof can be better understood if we provide the proof first only for the special case of the $C^*$-algebra of all $2 \times 2$ matrices and then turn to the proof of the general case.

\subsection{The case of the algebra of $2 \times 2$ matrices} \label{susec:2x2}

Let $I \subset \R$ be an open interval and $f$ be a function defined on $I$ that satisfies the conditions given in Theorem \ref{thm:main}.
Let $\cA$ be the $C^*$-algebra of all $2 \times 2$ complex matrices (which is denoted by $M_2(\C)$). Let $A \in M_2(\C)$ be a self-adjoint matrix with $\sigma(A)\subset I.$
\par
The proof of the direction \eqref{tul:cent}$\Longrightarrow$ \eqref{tul:loc-conv} is clear. If $A$ is central, that is, $A= \lambda I_2$ (where $I_2$ denotes the identity element of $M_2(\C)$) for some $\lambda \in I$, then $f(A) \leq \frac{1}{2}\ler{f(A+B)+f(A-B)}$ holds for every self-adjoint $B \in M_2(\C)$ (such that $\sigma(A+B) \cup \sigma(A-B) \subset I$) because of the convexity of $f$ as a scalar function.
\par
The interesting part is the proof of the direction \eqref{tul:loc-conv}$\Longrightarrow$\eqref{tul:cent}. We will prove it by contraposition, that is, we show that if $A$ is not central, then $f$ is not locally convex at the point $A.$ So assume that the self-adjoint matrix $A$ is not central, which means that it has two different eigenvalues, say, $x$ and $y$ in $I.$
\par
Let us use the formula for the (higher order) Fr\'echet derivatives of matrix valued functions defined by the functional calculus given by \emph{Hiai} and \emph{Petz} \cite[Thm. 3.33]{HP_book}. This formula is essentially based on the prior works of \emph{Daleckii} and \emph{Krein} \cite{dal-krein}, \emph{Bhatia} \cite{bhatia}, and \emph{Hiai} \cite{Hiai-m-anal}.

This formula gives us that if $A= \left[ \ba{cc} x & 0 \\ 0 & y \ea\right]$ and $B=\left[ \ba{cc} 1 & 1 \\ 1 & 1 \ea\right],$
then the second order Fr\'echet derivative of the function $f$ (defined by the functional calculus) at the point $A$ with arguments $(B,B)$ is
$$
\partial^2 f (A)(B,B)=2 \left[ \ba{cc} f^{[2]}[x,x,x]+f^{[2]}[x,x,y] & f^{[2]}[x,x,y]+f^{[2]}[x,y,y] \\ f^{[2]}[x,x,y]+f^{[2]}[x,y,y] & f^{[2]}[x,y,y]+f^{[2]}[y,y,y] \ea\right],
$$
where $f^{[2]}[\cdot,\cdot,\cdot]$ denotes the \emph{second divided difference} with respect to $f.$ (For the Fr\'echet derivatives, we use the notation of Hiai and Petz \cite{HP_book}.)
\par
It is well-known that
\begin{eqnarray}
& & \partial^2 f (A)(B,B)= \frac{\mathrm{d}^2}{\mathrm{d}t^2} f(A+tB)_{|t=0}\nonumber\\
& & =\lim_{t \to 0} \frac{1}{t^2}\ler{f(A+tB)-2f(A)+f(A-tB)}. \label{eq:2nd-der}
\end{eqnarray}
Now we show that $\partial^2 f (A)(B,B)$ is not positive semidefinite. Indeed, if $w=[1 \;\:-1]^\top$, then
\begin{eqnarray}
& & \inner{\partial^2 f (A)(B,B) w}{w}\nonumber\\
& & =2 \left(f^{[2]}[x,x,x]+f^{[2]}[x,x,y] - f^{[2]}[x,x,y]-f^{[2]}[x,y,y]\right.\nonumber\\
& & \qquad  \left. - f^{[2]}[x,x,y]-f^{[2]}[x,y,y] + f^{[2]}[x,y,y]+f^{[2]}[y,y,y]\right)\nonumber\\
& & = 2 \ler{f^{[2]}[x,x,x]-f^{[2]}[x,x,y]-f^{[2]}[x,y,y]+f^{[2]}[y,y,y]}, \label{eq:div-diff}
\end{eqnarray}
where $\inner{\cdot}{\cdot}$ denotes the inner product on $\C^2.$
Using the basic properties of the divided differences (which can be found e.g. in Section 3.4 of the book \cite{HP_book}) on can compute that the above expression \eqref{eq:div-diff} is equal to
\begin{eqnarray}
& & f''(x)-2 \frac{f'(x)-\frac{f(x)-f(y)}{x-y}}{x-y}-2 \frac{\frac{f(x)-f(y)}{x-y}-f'(y)}{x-y}+f''(y)\nonumber\\
& & =f''(x)-2 \frac{f'(x)-f'(y)}{x-y}+f''(y). \label{eq:egysz}
\end{eqnarray}
And the expression \eqref{eq:egysz} is negative by the strict concavity of the function $f''$ as one can see for example by the following calculation:
\begin{eqnarray}
& & 2 \ler{\frac{1}{2}\ler{f''(x)+f''(y)}-\frac{f'(x)-f'(y)}{x-y}}\nonumber\\
& & =2 \ler{\int_0^1 t f''(x)+(1-t)f''(y) \dt - \int_0^1 f''\ler{tx+(1-t)y}\dt}\nonumber\\
& & =2 \int_0^1 t f''(x)+(1-t)f''(y)- f''\ler{tx+(1-t)y}\dt.  \label{eq:int-repr}
\end{eqnarray}
The integrand in \eqref{eq:int-repr} is continuous in $t$ and is negative for every $0<t<1$ because $x \neq y$ and $f''$ is strictly concave, hence the above integral \eqref{eq:int-repr} is negative. So we deduced that $\inner{\partial^2 f (A)(B,B) w}{w}<0.$ (It is fair to remark that the above computation is essentially a possible proof of the well-known \emph{Hermite-Hadamard inequality}.)

\par
So, by \eqref{eq:2nd-der}, we have
$$
\inner{\lim_{t \to 0} \frac{1}{t^2}\ler{f(A+tB)-2f(A)+f(A-tB)} w}{w}<0.
$$
This means that
$$
\lim_{t \to 0} \frac{1}{t^2} \inner{\ler{f(A+tB)-2f(A)+f(A-tB)} w}{w}<0,
$$
so there exists some $t_0>0$ such that
\be \label{eq:delta-def}
\inner{\ler{f(A+t_0 B)-2f(A)+f(A-t_0 B)} w}{w}<0.
\ee
(For further use, let us denote the negative number in \eqref{eq:delta-def} by $-\delta.$)
So, we obtained that $f(A+t_0B)-2f(A)+f(A-t_0B)$ is not positive semidefinite, i.e.,
$$
0 \nleq f(A+t_0B)-2f(A)+f(A-t_0B),
$$
in other words,
$$
f(A) \nleq \frac{1}{2}\ler{f(A+t_0B)+f(A-t_0B)}.
$$
This means that $f$ is not locally convex at the point $A.$ The proof is done.

\subsection{The general case} \label{susec:general}
The proof of Theorem \ref{thm:main} in the case of a general $C^*$-algebra is heavily based on our arguments given in \cite{dv-mon}. For the convenience of the reader, we repeat some of the arguments of \cite{dv-mon} here in this subsection instead of referring to \cite{dv-mon} all the time.
\par

Also in this general case, the proof of the direction \eqref{tul:cent} $\Longrightarrow$ \eqref{tul:loc-conv} is easy. As $f$ is continuous and convex as a function of one real variable, the map $a \mapsto f(a)$ is also convex on any set of commuting self-adjoint elements of a $C^*$-algebra (provided that the expression $f(a)$ makes sense). So, centrality automatically implies local convexity.
\par
To prove the direction \eqref{tul:loc-conv} $\Longrightarrow$ \eqref{tul:cent}, we use contraposition again. Assume that $a \in \cA_s, \, \sigma(a) \subset I$ and $a$ is not central, that is, $aa'-a'a \neq 0$ for some $a' \in \cA.$ Then, by \cite[10.2.4. Corollary]{kad-ring-2}, there exists an irreducible representation $\pi: \cA \rightarrow \cB(\cH)$ such that $\pi\ler{aa'-a'a} \neq 0,$ that is, $\pi(a) \pi\ler{a'}\neq \pi\ler{a'}\pi(a).$ Let us fix this irreducible representation $\pi.$ So, $\pi(a)$ is a non-central self-adjoint (and hence normal) element of $\cB(\cH)$ with $\sigma\ler{\pi(a)}\subset I$ (as a representation does not increase the spectrum). By the non-centrality, $\sigma\ler{\pi(a)}$ has at least two elements, and by the normality, every element of $\sigma\ler{\pi(a)}$ is an approximate eigenvalue \cite[3.2.13. Lemma]{kad-ring-1}. Let $x$ and $y$ be two different elements of $\sigma\ler{\pi(a)},$ and let $\{u_n\}_{n \in \N} \subset \cH$ and $\{v_n\}_{n \in \N} \subset \cH$ satisfy
\begin{eqnarray*}
& & \lim_{n \to \infty} (\pi(a) u_n - x u_n) =0, \,
\lim_{n \to \infty} (\pi(a) v_n - y v_n) =0,\\
& & \inner{u_m}{v_n}=0 \;\:\mbox{for all}\;\: m,n \in \N.
\end{eqnarray*}
(As $x\neq y,$ the approximate eigenvetors can be chosen to be orthogonal.) Set $\cK_n:=\mathrm{span}\{u_n, v_n\}$ and let $E_n$ be the orthoprojection onto the closed subspace $\cK_n^{\perp} \subset \cH.$ Let
$$
\psi_n(a):=x u_n \otimes u_n + y v_n \otimes v_n + E_n \pi(a)E_n.
$$
We intend to show that
$$ 
\lim_{n \to \infty} \psi_n(a)=\pi(a)
$$
in the operator norm topology.
Let $h$ be an arbitrary non-zero element of $\cH$ and consider the orthogonal decompositions $h=h_1^{(n)}+h_2^{(n)},$ where $h_1^{(n)} \in \cK_n$ and $h_2^{(n)} \in \cK_n^{\perp}$ for any $n \in \N.$
Let us introduce the symbols $\eps_{u,n}:=\pi(a) u_n-x u_n$ and $\eps_{v,n}:=\pi(a) v_n-y v_n$ and recall that $\lim_{n \to \infty} \eps_{u,n}=0$ and $\lim_{n \to \infty} \eps_{v,n}=0$ in the standard topology of the Hilbert space $\cH.$
Now,
\begin{eqnarray*}
& & \frac{1}{\norm{h}}
\norm{\ler{\pi(a)-\psi_n(a)}h}\\
& & \leq
\frac{1}{\norm{h}}
\norm{\ler{\pi(a)-\psi_n(a)}h_1^{(n)}}
+\frac{1}{\norm{h}}\norm{\ler{\pi(a)-\psi_n(a)}h_2^{(n)}}.
\end{eqnarray*}
Both the first and the second term of the right hand side of the above inequality are bounded by the term $\norm{\eps_{u,n}}+\norm{\eps_{v,n}}$ because
\begin{eqnarray*}
& & \frac{1}{\norm{h}}\norm{\ler{\pi(a)-\psi_n(a)}h_1^{(n)}}
 =\frac{1}{\norm{h}}\norm{\ler{\pi(a)-\psi_n(a)}\ler{\alpha_n u_n +\beta_n v_n}}\\
& &
 =\frac{1}{\norm{h}}
\norm{\alpha_n x u_n + \alpha_n \eps_{u,n}-x \alpha_n u_n + \beta_n y v_n + \beta_n \eps_{v,n}-y \beta_n v_n}\\
& &  \leq \frac{\abs{\alpha_n}}{\norm{h}}
\norm{\eps_{u,n}}+\frac{\abs{\beta_n}}{\norm{h}}\norm{\eps_{v,n}}
\leq \norm{\eps_{u,n}}+\norm{\eps_{v,n}}
\end{eqnarray*}
as the sequences $\left\{\abs{\alpha_n}\right\}$ and $\left\{\abs{\beta_n}\right\}$ are obviously bounded by $\norm{h},$
and
\begin{eqnarray*}
& & \frac{1}{\norm{h}}
\norm{\ler{\pi(a)-\psi_n(a)}h_2^{(n)}}=
\frac{1}{\norm{h}}
\norm{\ler{I_\cH-E_n}\pi(a)h_2^{(n)}}\\
& & =\frac{1}{\norm{h}}
\norm{\ler{u_n \otimes u_n + v_n \otimes v_n}\pi(a)h_2^{(n)}}\\
& & =\frac{1}{\norm{h}}
\norm{\inner{\pi(a)h_2^{(n)}}{u_n}u_n+\inner{\pi(a)h_2^{(n)}}{v_n}v_n}\\
& & =
\frac{1}{\norm{h}}
\norm{\inner{h_2^{(n)}}{\pi(a)u_n}u_n+\inner{h_2^{(n)}}{\pi(a) v_n}v_n}\\
& & \leq
\frac{1}{\norm{h}}
\abs{\inner{h_2^{(n)}}{x u_n +\eps_{u,n}}}+ \frac{1}{\norm{h}}
\abs{\inner{h_2^{(n)}}{y v_n +\eps_{v,n}}}\\
& & =
\frac{1}{\norm{h}}
\abs{\inner{h_2^{(n)}}{\eps_{u,n}}}
+\frac{1}{\norm{h}}
\abs{\inner{h_2^{(n)}}{\eps_{v,n}}}\\
& & \leq
\frac{\norm{h_2^{(n)}}}{\norm{h}}\norm{\eps_{u,n}}
+\frac{\norm{h_2^{(n)}}}{\norm{h}}\norm{\eps_{v,n}}
\leq \norm{\eps_{u,n}}+\norm{\eps_{v,n}}.
\end{eqnarray*}
We used that $a$ is self-adjoint and that hence so is $\pi(a).$
So, we found that
$$
\mathrm{sup}
\left\{
\frac{1}{\norm{h}}
\norm{\ler{\pi(a)-\psi_n(a)}h}
\middle|
h \in \cH\setminus \{0\} \right\}
\leq 2 \ler{\norm{\eps_{u,n}}+\norm{\eps_{v,n}}} \to 0,
$$
which means that $\psi_n(a)$ tends to $\pi(a)$ in the operator norm topology.

\par
Let us use the notation $B_n:=(u_n+v_n) \otimes (u_n+v_n)$ and $w_n:= u_n - v_n.$ By the result of Subsection \ref{susec:2x2} (the proof for the case of $\cA=M_2(\C)$) we have
\begin{eqnarray}
& & \inner{
\ler{f\ler{\psi_n(a)+t_0 B_n}-2 f\ler{\psi_n(a)}+f\ler{\psi_n(a)-t_0 B_n}} w_n}{w_n}\nonumber\\
& & =-\delta<0, \label{deltas}
\end{eqnarray}
where $t_0$ is the same as in \eqref{eq:delta-def} and $-\delta$ is the left hand side of \eqref{eq:delta-def}, for any $n \in \N.$ That is, the left hand side of \eqref{deltas} is independent of $n.$
\par
The operator $B_n$ is a self-adjoint element of $\cB(\cH)$ and $\cK_n$ is a finite dimensional subspace of $\cH,$ hence by Kadison's transitivity theorem \cite[10.2.1. Theorem]{kad-ring-2}, there exists a self-adjoint $b_n \in \cA$ such that
$$
\pi\ler{b_n}_{|\cK_n}=B_{n |\cK_n}.
$$
So, we can rewrite \eqref{deltas} as
\begin{eqnarray}
& & \inner{
\ler{f\ler{\psi_n(a)+t_0 \pi\ler{b_n}}-2 f\ler{\psi_n(a)}+f\ler{\psi_n(a)-t_0 \pi\ler{b_n}}} w_n}{w_n}\nonumber\\
& & =-\delta<0, \label{uj_delt}
\end{eqnarray}
A standard continuity argument which is based on the fact that $\psi_n(a)$ tends to $\pi(a)$ in the operator norm topology shows that
\be \label{fkonv1}
\lim_{n \to \infty} \norm{f\ler{\psi_n(a)}-f\ler{\pi(a)}}=0.
\ee
Moreover, by Kadison's transitivity theorem, the sequence $\pi \ler{b_n}$ is bounded (for details, the reader should consult the proof of \cite[5.4.3. Theorem]{kad-ring-1}), and hence
\be \label{fkonv2}
\lim_{n \to \infty} \norm{f\ler{\psi_n(a) \pm t_0\pi \ler{b_n}}-f\ler{\pi(a) \pm t_0\pi \ler{b_n}}}
=0
\ee
also holds.
By \eqref{fkonv1} and \eqref{fkonv2}, for any $\delta>0$ one can find $n_0 \in \N$ such that for $n>n_0$ we have
$$
\norm{f\ler{\psi_n(a)}-f\ler{\pi(a)}}<\frac{1}{16} \delta
$$
and
$$
\norm{f\ler{\psi_n(a) \pm t_0\pi \ler{b_n}}-f\ler{\pi\ler{a \pm t_0 b_n}}} <\frac{1}{16} \delta.
$$
Therefore, by \eqref{uj_delt}, for $n>n_0,$ the inequality
$$ 
\inner{
\ler{f\ler{\pi\ler{a+t_0 b_n}}-2 f\ler{\pi(a)}+f\ler{\pi\ler{a-t_0 b_n}}} w_n}{w_n}<-\delta/2<0,
$$
holds.
In other words,
$$
f\ler{\pi(a)} \nleq \frac{1}{2} \ler{f\ler{\pi\ler{a+t_0 b_n}}+f\ler{\pi\ler{a-t_0 b_n}}}
$$
or equivalently (as functional calculus commutes with every representation of a $C^*$-algebra),
$$
\pi\ler{f(a)} \nleq \pi \ler{\frac{1}{2} \ler{f\ler{a+t_0 b_n}+f\ler{a-t_0 b_n}}}.
$$
Any representation of a $C^*$-algebra preserves the semidefinite order, hence this means that
$$
f(a) \nleq \frac{1}{2} \ler{f\ler{a+t_0 b_n}+f\ler{a-t_0 b_n}},
$$
which means that $f$ is not locally convex at $a.$ The proof is done.

\subsection*{Acknowledgement}
The author is grateful to Lajos Moln\'ar for proposing the problem discussed in this paper and for great conversations about this topic and about earlier versions of this paper. The author is grateful to Albrecht B\"ottcher for suggestions that helped to improve the presentation of this paper.

\end{document}